\documentclass[12pt]{amsproc}
\usepackage{amscd}
\usepackage{amssymb}

\newtheorem{theorem}{Theorem}[section]

\newtheorem{lemma}[theorem]{Lemma}
\newtheorem{proposition}[theorem]{Proposition}

\theoremstyle{definition}
\newtheorem{definition}[theorem]{Definition}

\newtheorem{examples}[theorem]{Examples}
\newtheorem{ass}[theorem]{Assumption}

\newtheorem{rem}[theorem]{Remark}

\newtheorem{prob}[theorem]{Problem}

\theoremstyle{remark}

\numberwithin{equation}{section}

\newsavebox{\SmallMathBox}

\DeclareRobustCommand*{\nicefrac}[2]{\ifmmode\mathnicefrac{#1}
{ #2}%
  \else\textnicefrac{#1}{#2}\fi}
\newcommand*{\textnicefrac}[2]{\check@mathfonts%
\mbox{\raisebox{.5ex}{\fontsize\sf@size\z@\selectfont#1}\kern-.
1em%
/\kern-.1em\raisebox{- .25ex}{\fontsize\sf@size\z@\selectfont#2}
}}
\newcommand*{\mathnicefrac}[2]{%
  \mathchoice
    {\m@fr@c{\scriptstyle}{#1}{#2}}
    {\m@fr@c{\scriptstyle}{#1}{#2}}
    {\m@fr@c{\scriptscriptstyle}{#1}{#2}}
    {\m@fr@c{\scriptscriptstyle}{#1}{#2}}}

\newcommand{\abs}[1]{\lvert#1\rvert}

\def\Ci{C^\infty}

\def\lla{\langle}

\def\noi{\noindent}
\newcommand{\norm}[1]{\lVert#1\rVert}

\def\ol{\overline}

\def\rra{\rangle}
\def\sqm1{\sqrt{-1}}

\def\tand{\mbox{\ \rm  and }}

\def\too{\longrightarrow}

\def\wt{\widetilde}

\def\={\cong}
\def\>{\supset}
\def\<{\subset}
\def\ii{^{-1}}
\def\12{\frac{1}{2}}
\def\0{^{\circ}}

\def\CC{{\mathbb C}}

\def\RR{{\mathbb R}}

\def\Aa{{\mathcal A}}
\def\Bb{{\mathcal B}}
\def\Cc{{\mathcal C}}

\def\Ff{{\mathcal F}}
\def\Gg{{\mathcal G}}

\def\Ll{{\mathcal L}}

\def\Nn{{\mathcal N}}

\def\Uu{{\mathcal U}}

\def\GGG{{\mathfrak G}}

\def\C{\CC}

\def\f{\varphi}

\def\g{\gamma}

\def\la{\lambda}

\def\R{\RR}

\def\Si{\Sigma}

\def\w{\omega}

\DeclareMathOperator{\dom}{dom}

\DeclareMathOperator{\id}{id}

\def\index{\mbox{\rm index\,}}

\DeclareMathOperator{\Mas}{Mas} \DeclareMathOperator{\mmax}{max}

\DeclareMathOperator{\ran}{im} \DeclareMathOperator{\range}{im}

\DeclareMathOperator{\sa}{sa}
\DeclareMathOperator{\SF}{sf}
\DeclareMathOperator{\sign}{sign}

\begin{document}
\setcounter{page}{1} \setcounter{tocdepth}{2}

\renewcommand*{\labelenumi}{%
   (\roman{enumi})}

\title
{General Spectral Flow Formula for Fixed Maximal Domain}

\thanks{This work was supported in part by The Danish Science Research Council,
SNF grant 21-02-0446. The second author is partially supported by
FANEDD 200215, 973 Program of MOST, Fok Ying Tung Edu. Funds
91002, LPMC of MOE of China, and Nankai University.}

\author{Bernhelm Booss-Bavnbek}
\address{Institut for matematik og fysik\\ Roskilde
University, 4000 Ros\-kilde, Denmark} \email{booss@mmf.ruc.dk}
\urladdr{http://imfufa.ruc.dk/$\sim$Booss}

\author{Chaofeng Zhu}
\address{Nankai Institute of Mathematics, Key Lab of Pure Mathematics
and Combinatorics of Ministry of Education, Nankai University,
Tianjin 300071, People's Republic of China}
\email{zhucf@nankai.edu.cn}

\date{}

\begin{abstract}
We consider a continuous curve of linear elliptic formally
self-adjoint differential operators of first order with smooth
coefficients over a compact Riemannian manifold with  boundary
together with a continuous curve of global elliptic boundary value
problems. We express the spectral flow of the resulting continuous
family of (unbounded) self-adjoint Fredholm operators in terms of
the Maslov index of two related curves of Lagrangian spaces. One
curve is given by the varying domains, the other by the Cauchy
data spaces. We provide rigorous definitions of the underlying
concepts of spectral theory and symplectic analysis and give a
full (and surprisingly short) proof of our General Spectral Flow
Formula for the case of fixed maximal domain. As a side result, we
establish local stability of weak inner unique continuation
property (UCP) and explain its role for parameter dependent
spectral theory.
\end{abstract}

\maketitle




\section{Statement of the problem and main result}\label{ss:statement}

\subsection{Statement of the problem}
Roughly speaking, the spectral flow counts the net number of
eigenvalues changing from the negative real half axis to the
non-negative one. The definition goes back to a famous paper by
M.~Atiyah, V.~Patodi, and I.~Singer \cite{AtPaSi75}, and was made
rigorous by J.~Phillips \cite{Ph96} for continuous paths of
bounded self-adjoint Fredholm operators, by K.P. Wojciechowski
\cite{Wo85} and C. Zhu and Y. Long \cite{ZhLo99} in various
non-self-adjoint cases, and by B. Booss-Bavnbek, M. Lesch, and J.
Phillips \cite{BoLePh01} in the unbounded self-adjoint case. We
shall give a rigorous definition of spectral flow, most suitable
for our purpose, below in Subsection \ref{ss:spectral-flow}
together with a review of its basic properties. For a definition of spectral flow
admitting zero in the continuous spectrum, we refer to A. Carey and J. Phillips \cite{CaPh03}.

\smallskip

In various branches of mathematics one is interested in the
calculation of the spectral flow of a continuous family of closed
densely defined (not necessarily bounded) self-adjoint Fredholm
operators in a fixed Hilbert space. We consider the following
typical problem of this kind.

\begin{ass}\label{a:basic-assumptions}
Let $\{A_s:\Ci(M;E)\to\Ci(M;E)\}_{s\in [0,1]}$ be a family of
formally self-adjoint linear elliptic differential operators of
first order with continuously varying smooth coefficients over a smooth
compact Riemannian manifold $M$ with boundary $\Si$, acting on
sections of a Hermitian vector bundle $E$ over $M$. Let $\{P_s\}$
be a continuous family of orthogonal pseudodifferential
projections in $L^2(\Si;E|_{\Si})$. Define $A_{s,P_s}$ to be the
unbounded operator in $L^2(M;E)$ with domain
\begin{equation}\label{e:domain}
D_s:=\{x\in H^1(M;E) \mid P_s(\gamma(x))=0\},
\end{equation}
where
\begin{equation}\label{e:gamma-sobolev}
\gamma:H^1(M;E)\to H^{\frac{1}{2}}(\Si;E|_{\Si})
\end{equation}
denotes the (continuous) trace map from the first Sobolev space over the whole
manifold to the $\12$ Sobolev space over the boundary.
(Note that in this paper the symbols $x$ and $y$ do not denote points of the underlying
manifolds $M$ or $\Si$, but points in Hilbert spaces, sections of vector bundles, etc.,
following the conventions of functional analysis and dynamical systems.)
Assume that
each $P_s$ defines a self-adjoint elliptic boundary condition for
$A_s$, i.e., $A_{s,P_s}$ is a self-adjoint Fredholm operator for
each $s\in[0,1]$.
\end{ass}

Then the spectral flow $\SF\{A_{s,P_s}\ ; \ s\in [0,1]\}$ or,
shortly, $\SF\{A_{s,P_s} \}$ is well defined. As a spectral
invariant it is essentially a {\em quantum} variable which one may
not always be able to determine directly by eigenvalue
calculations. As an alternative, one is looking for a {\em
classical} method of calculating the spectral flow.
There are two different approaches. One setting expresses the spectral
flow (of a loop of Dirac operators on a closed manifold) as an integral
over a 1-form induced by the heat kernel (for a review see \cite{CaPh03}).
The other setting is reduction to the boundary, i.e., one expresses the
spectral flow (of a path of self-adjoint boundary value problems on a
compact manifold with boundary) in terms of
the intersection geometry of the solution spaces of the
homogeneous differential equations and the boundary conditions.
That is the approach we shall follow in this paper.

\begin{prob} \label{pr:sf}Give a {\em classical} method of calculating the spectral
flow of the family $\{A_{s,P_s}\}$ by reduction to the
boundary, i.e., a method not involving
the determination of the spectrum near $0$ and yielding an expression on $\Sigma$.
\end{prob}

The preceding spectral flow calculation problem is formulated for
families by analogy with {\em Bojarski's Theorem} for single
operators which expresses the index (which is the difference
between the multiplicities of the 0-eigenvalue of the original and
the formally adjoint problem and so a priori a quantum or spectral
invariant) of an elliptic operator over a closed partitioned
manifold $M=M_-\cup_{\Sigma} M_+$ by the index of the Fredholm
pair of Cauchy data spaces from two sides along the hypersurface
$\Si$ (which are classical objects, see Bojarski \cite{Bo79} and
Booss and Wojciechowski \cite[Chapter 24]{BoWo93}).

\smallskip

\subsection{General functional analytic setting and announcement of the General
Spectral Flow Formula} Now we translate our problem into a
functional analytic setting. For any such family there are three
geometrically defined relevant Hilbert spaces of global sections
which remain fixed under variation of the coefficients of the
operators and under variation of the boundary conditions:
\begin{equation}\label{e:fa-setting}
L^2(M;E), \quad H_0^1(M;E), \quad {\rm and}\;\;H^1(M;E).
\end{equation}
Here $H_0^1(M;E)$ denotes the closure of $\Ci_0(M\setminus
\Si;E)$ in the first Sobolev space $H^1(M;E)$, where
$\Ci_0(M\setminus \Si;E)$ denotes the smooth sections with support
in the interior of $M\setminus \Si$. Since the trace map
$\gamma:H^1(M;E)\to H^{\frac{1}{2}}(\Si;E|_{\Si})$ is continuous,
we have $H_0^1(M;E)=\ker\g$, i.e., the space $H^1_0(M;E)$ consists
exactly of the elements of $H^1(M;E)$ which vanish on the boundary
$\Si$\,.

For each $s\in[0,1]$, we shall denote the unbounded operator $A_s$
acting in $L^2(M;E)$ with domain $H_0^1(M;E)$ also by $A_s$\,.
Since the differential operator $A_s$ is elliptic, the unbounded
operator $A_s$ is closed by {\em G{\aa}rding's inequality}
\begin{equation}\label{e:gaarding}
\norm{x}_{H^1(M;E)} \le
C\bigl(\norm{x}_{L^2(M;E)}+\norm{A_sx}_{L^2(M;E)}\bigr)
\quad\text{ for $x\in H_0^1(M;E)$}\,.
\end{equation}

Denote by $\dom(A)$ the domain of an operator $A$, by $A^*$ the
adjoint operator of $A$, and
\begin{equation}\label{e:dom-max}
D_{\mmax}(A):=\dom(A^*).
\end{equation}
Since $A$ is closed and
symmetric, it follows that $D_{\mmax}(A) = \{x\in L^2(M;E) \mid Ax
\in L^2(M;E)\}$ with $Ax$ taken in the distributional sense. For
$A_s$ formally self-adjoint, it follows immediately that
$H^1(M;E)\subset D_{\mmax}(A_s)$ and that $A_s$
(with domain $H^1_0(M;E)$) is symmetric.

\medskip

In local coordinates, we view each coefficient of $A_s$ as a map which maps
$s\in[0,1]$ to a continuous section (which is actually smooth). Then the
continuity of the curve $\{A_s\}_{s\in[0,1]}$ in the sense of {\em
continuously varying coefficients} implies the continuity of the
curve
\begin{equation}\label{e:bounded-operator-continuity}
[0,1]\ni s\mapsto A_s^*|_{H^1(M;E)}\in \Bb(H^1(M;E),L^2(M;E))\,,
\end{equation}
as a curve of {\em bounded operators} from $H^1(M;E)$ to
$L^2(M;E)$.

We denote by $Q_s:L^2(\Si;E|_{\Si})\to L^2(\Si;E|_{\Si})$ the {\em
Calder{\'o}n projection}. It is a projection onto the {\em Cauchy
data space} of $A_s^*$ which is defined as the $L^2$-closure of
$\g(\ker (A_s^*|_{H^1(M;E)})$. It can be described as a
pseudodifferential operator, e.g., when continuing $A_s$ to an
elliptic operator on a closed manifold $\wt M\> M$, see R.T.
Seeley \cite[Sections 4 and 8]{Se66} and \cite[Chapter
12]{BoWo93}. For an alternative canonical construction based on a
natural boundary value problem and avoiding the choices of closing
the manifold and continuing the operator, see B. Himpel, P. Kirk,
and M. Lesch \cite[Section 3]{HiKiLe04} and recent joint work of
the authors with M. Lesch \cite{BoLeZh05}.

\smallskip

For each $s\in[0,1]$, there is a natural (strong) symplectic form
$\omega_s$ on the quotient space $D_{\mmax}(A_s)/H_0^1(M;E)$
induced by {\em Green's form} of $A_s$ as
\begin{equation}\label{e:green-form}
\omega_s(\gamma(x),\gamma(y)):=\lla A_s^*x,y\rra -\lla
x,A_s^*y\rra\,, \quad x,y\in D_{\mmax}(A_s).
\end{equation}
Here $\g$ denotes the natural projection
\[
D_{\mmax}(A_s) \to D_{\mmax}(A_s)/H^1_0(M;E)\,.
\]
Identifying the quotient space
$D_{\mmax}(A_s)/H^1_0(M;E)$ with a subspace of the Sobolev (distribution)
space $H^{-1/2}(\Si;E|_{\Si})$, we obtain that this $\g$ extends
the Sobolev trace map of \eqref{e:gamma-sobolev}. A rigorous
definition of {\em symplectic structures} and {\em Lagrangian
subspaces} will be given below in Subsection \ref{ss:symplectic}.

For our formally self-adjoint differential operators of first
order, we have an explicit description of the form in
\eqref{e:green-form}, restricted to $H^1(M;E)$\,, by Stokes'
Theorem
\begin{equation}\label{e:green-form-explicit}
\omega_s(\gamma(x),\gamma(y)) =-\int_{\Si}
\lla\sigma_1(A_s)(\cdot,dt) \bigl(x|_{\Si}\bigr),y|_{\Si}\rra\/
d{\rm vol}_{\Si}\,,
\end{equation}
where $\sigma_1(A_s)(\cdot,dt )$ denotes the principal symbol of
$A_s$ at the boundary, taken in inner (co-)normal direction $dt$.
Notice that we do not require that the manifold $M$ is orientable:
for our application of Stokes' Theorem it suffices that any collar
neighborhood of $\Si$ in $M$ is oriented by the normal structure.
Then the form $\omega_s|_{H^1(M;E)}$ of
\eqref{e:green-form-explicit} extends to a (strong) symplectic
structure $\ol{\omega}_s$ on $L^2(\Si;E|_{\Si})$. One can show
that $\omega_s|_{H^1(M;E)/H_0^1(M;E)}$ is a {\em weak} (but not
strong) symplectic form on the Hilbert space $H^1(M;E)/H_0^1(M;E)
\cong H^{\12}(\Si;E|_{\Si})$ (cf. Booss and Zhu \cite[Remark
1.6b]{BoZh04}).

We have $H^1(M;E)=D_{\mmax}(A_s)$ if and only if $\dim M=1$. For
higher dimensional case, the strict inclusion $H^1(M;E)\subset
D_{\mmax}(A_s)$ and the weakness of $\omega_s|_{H^1(M;E)}$ causes
technical difficulties.


However, we still have the following theorem (cf. Theorem 0.1 of
\cite{BoZh04}).

\begin{theorem}[General Spectral Flow Formula] \label{t:gsff}
Let $\{A_s\}_{s\in[0,1]}$ and $\{P_s\}_{s\in[0,1]}$ be operator
families like in Assumption \ref{a:basic-assumptions}. We assume
that $\{\ker P_s\}_{s\in[0,1]}$ is a continuous family of
Lagrangian subspaces in $(H,\ol{\omega}_s)$. If $A_s$ satisfies
weak inner UCP, i.e., $\ker A_s=\{0\}$ for each $s\in [0,1]$, we
have:

\noi {\rm (a)} The family $\{A_{s,P_s}\}_{s\in[0,1]}$ of closed
self-adjoint Fredholm operators on $X$ is a continuous family
$($in the gap norm, or equivalently, in the projection norm$)$.

\noi {\rm (b)} The Cauchy data spaces $\range Q_s$ are Lagrangian
subspaces in the weak symplectic Hilbert space $(
H^{\12}(\Si;E|_{\Si}),\ol{\omega}_s)$ and form a continuous family
in $H^{\12}(\Si;E|_{\Si})$ for $s\in[0,1]$.

\noi {\rm (c)} Finally, the following formula holds:
\begin{equation}\label{e:gsff}
\SF\{A_{s,P_s}\}=-\Mas\{\ker P_s,\range Q_s\},
\end{equation}
where the spectral flow $\SF$ and the Maslov index $\Mas$ are
defined by Definitions {\rm \ref{d:sf-sa}} and {\rm
\ref{d:maslov}} below respectively.
\end{theorem}

\begin{rem}\label{r:gsff}
(a) The General Spectral Flow Formula contains and generalizes
all previously known spectral flow formulae, as given by M. Morse
\cite{Mo34}, W. Ambrose \cite{Am61}, J.J. Duistermaat \cite{Du76},
A. Floer \cite{Fl88}, P. Piccione and D.V. Tausk \cite{PiTa00} and
\cite{PiTa02}, and C. Zhu \cite{Zh00} and \cite{Zh01} for the
1-dimensional setting of the study of geodesics, and for the
higher dimensional setting the formulae given by T. Yoshida
\cite{Yo91}, L. Nicolaescu \cite{Ni95}, S.E. Cappell, R. Lee, and
E.Y. Miller \cite{CaLeMi96}, B. Booss, K. Furutani, and N. Otsuki
\cite{BoFu98} and \cite{BoFuOt01}, and P. Kirk and M. Lesch
\cite{KiLe00}.

(b) The main difference to \cite{BoFu98} and \cite{BoFuOt01}
is that we admit varying maximal domain and varying Fredholm
domain. The main difference to \cite{KiLe00} is that we admit more
general operators than Dirac type operators with constant
coefficients in normal direction close to the boundary.

(c) The proof of the above theorem is rather technical and
complicated. In this review article, we only prove the following
fixed maximal domain case which completely covers all above cited
one-dimensional cases (cf. Corollary 2.14 in \cite{BoZh04}).
Moreover, it contains \cite{BoFu98} and \cite{BoFuOt01} and
generalizes it to varying Fredholm domains, and contains
\cite{KiLe00} for the case of fixed maximal domain and generalizes
it under that restriction to more general operator families.
\end{rem}

\smallskip

\subsection{Statement of the result for fixed maximal
domain}\label{ss:statement-fixd} Let $X$ be a Hilbert space, and
$D_m\subset D_{\mmax}$ be two dense linear subspaces of $X$. Let
$\{A_s\}_{s\in[0,1]}$ be a family of symmetric densely defined
operators in $\Cc(X)$ with domain $\dom(A_s)=D_m$. Here we denote
by $\Cc(X)$ all closed operators in $X$. Assume that
$\dom(A_s^*)=D_{\mmax}$\,, i.e., the domain of the maximal
symmetric extension $A^*_s$ of $A_s$ is independent of $s$.

We recall from \cite{BoFu98} (see also B. Lawruk, J.
{\'S}niatycki, and W.M. Tulczyjew \cite{LaSnTu75} for early
investigation of symplectic structures and boundary value
problems) for each $s\in[0,1]$:

\begin{enumerate}

\item The space $D_{\mmax}$ is a Hilbert space with the graph inner
product
\begin{equation}\label{e:graph_inner_product1}
\lla x,y\rra_{\GGG_s} := \lla x,y\rra_X + \lla A_s^*x,A_s^*y\rra_X
\quad\text{ for $x,y\in D_{\mmax}$}\,.
\end{equation}

\item The space $D_m$ is a closed subspace in the graph norm and
the quotient space $D_{\mmax}/D_m$ is a strong symplectic Hilbert
space with the (bounded) symplectic form induced by Green's form
\begin{equation}\label{e:symplectic_green1}
\omega_s(x+D_m,y+D_m) := \lla A_s^*x,y\rra_X - \lla x,A_s^*y\rra_X
\quad\text{ for $x,y\in D_{\mmax}$}\,.
\end{equation}

\item If $A_s$ admits a self-adjoint Fredholm extension
$A_{s,D_s}:=A_s^*|_{D_s}$ with domain $D_s$, then the {\em
natural Cauchy data space} $(\ker A_s^*+D_m)/D_m$ is a Lagrangian
subspace of $(D_{\mmax}/D_m,\omega_s)$\,.

\item Moreover, self-adjoint Fredholm extensions are characterized
by the property of the domain $D_s$ that $(D_s+D_m)/D_m$ is a
Lagrangian subspace of $(D_{\mmax}/D_m,\omega_s)$ and forms a
Fredholm pair with $(\ker A_s^*+D_m)/D_m$\,.

\item We denote the natural projection (which is independent of
$s$) by
\[
\g:D_{\mmax}\too D_{\mmax}/D_m\/.
\]
\end{enumerate}

The main result of this paper is the following theorem which
reproves parts of the preceding list.

\begin{theorem}[General Spectral Flow Formula for fixed maximal domain]
\label{t:fixm} We assume that on $D_{\mmax}$ the graph norms
induced by $A_s$, $0\le s\le 1$ are mutually equivalent. Then we
fix a graph norm $\Gg$ on $D_{\mmax}$ induced by $A_0$. Assume
that $\{A_s^*:D_{\mmax}\to X\}$ is a continuous family of bounded
operators and each $A_s$ is injective. Let $\{D_s/D_m\}$ be a
continuous family of Lagrangian subspaces of
$(D_{\mmax}/D_m,\w_s)$, such that each $A_{s,D_s}$ is a Fredholm
operator. Then:

\noi {\rm (a)} Each $\bigl(D_s/D_m,\gamma(\ker(A_s^*)\bigr)$ is a
Fredholm pair in $D_{\mmax}/D_m$\/.

\noi {\rm (b)} Each Cauchy data space $\g(\ker A_s^*)$ is a
Lagrangian subspace of $(D_{\mmax}/D_m,\omega_s)$\,.

\noi {\rm (c)} The family $\{\g(\ker A_s^*)\}$ is a continuous
family in $D_{\mmax}/D_m$\,.

\noi {\rm (d)} The family $\bigl\{A_{s,D_s}\bigr\}$ is a
continuous family of self-adjoint Fredholm operators in $\Cc(X)$.

\noi {\rm (e)} Finally, we have
\begin{equation}\label{e:fixd}
\SF\{A_{s,D_s}\} = -\Mas \{\g(D_s),\g(\ker A_s^*)\}.
\end{equation}
\end{theorem}

\smallskip

\noi{\em Acknowledgement.} The first author thanks the organizers
Jan Kubarski, Tomasz Rybicki, and Robert Wolak of the 6th
Conference on
 Geometry and Topology of Manifolds (Krynica, Poland, May 2-8, 2004) for the opportunity
to present the main ideas and various ramifications of this paper
in a mini-course of four hours.
We both thank the referee for corrections,
thoughtful comments, and helpful suggestions which led to many improvements.
The referee clearly went beyond the call of duty, and we are indebted.

\medskip

\section{Definition of spectral flow and Maslov index}

\subsection{Spectral flow, revisited and generalized}\label{ss:spectral-flow}

Let $X$ be a Hilbert space. For a self-adjoint Fredholm operator
$A\in \Cc(X)$, there exists a unique orthogonal decomposition
\begin{equation}\label{e:spectral-decomposition}
X=X^+(A)\oplus X^0(A)\oplus X^-(A)
\end{equation}
such that $X^+(A)$, $X^0(A)$ and $X^-(A)$ are invariant subspaces
associated to $A$, and $A|_{X^+(A)}$, $A|_{X^0(A)}$ and
$A|_{X^-(A)}$ are positive definite, zero and negative definite
respectively. We introduce vanishing, natural, or infinite numbers
$$
m^+(A):=\dim X^+(A),\ m^0(A):=\dim X^0(A),\ m^-(A):=\dim X^-(A),
$$
and call them {\em Morse positive index}, {\em nullity} and {\em
Morse index} of $A$ respectively.
 For finite-dimensional $X$, the {\em signature} of $A$ is defined by
$\sign(A)=m^+(A)-m^-(A)$ which yields an integer. The {\em APS
projection} $Q_A$ (where APS stands for Atiyah-Patodi-Singer) is
defined by
$$Q_A(x^++x^0+x^-):=x^++x^0,$$
for all $x^+\in X^+(A), x_0\in X^0(A), x^-\in X^-(A)$.

\medskip

Let $\{A_s\}$, $0\le s\le 1$ be a continuous family of
self-adjoint Fredholm operators. The spectral flow $\SF\{A_s\}$ of
the family should be equal to $m^-(A_0)-m^-(A_1)$ if $\dim
X<+\infty$. We will generalize this definition to general $X$.

For each $t\in[0,1]$, there exists a bounded open neighborhood
$N_t$ of $0$ such that $\partial N_t$ is of class $C^1$,
$\sigma(A_t)\cap\partial N_t=\emptyset$, and $P(A_t,N_t)$ is a
finite rank projection. Here we denote the spectrum of a closed
operator $A$ by $\sigma(A)$, and the spectral projection by
$$P(A,N):=-\frac{1}{2\pi\sqrt{-1}}\int_{\partial N}(A-zI)^{-1}dz$$
if $N$ is a bounded open subset of $\C$ with $C^1$ boundary and
$\partial N\cap \sigma(A)=\emptyset$. The orientation of $\partial
N$ is chosen to make $N$ stay on the left side of $\partial N$.
Since the family $\{A_s\}$, $0\le s\le 1$ is continuous, there
exists a $\delta(t)>0$ for each $t\in[0,1]$ such that
$$\sigma(A_s)\cap \partial N_t=\emptyset,\;{\rm for}\;{\rm
all}\;s\in(t-\delta(t),t+\delta(t))\cap[0,1].$$ Then
\[
\bigl\{P(A_s,N_t)\bigr\}_{s\in (t-\delta(t),t+\delta(t))\cap[0,1]}
\quad \text{for fixed $t\in[0,1]$},
\]
is a continuous family of orthogonal projections. By Lemma I.4.10
in Kato \cite{Ka80}, they have the same rank. We denote by
$A(s,t)$ the operator $A_s$ acting on the finite-dimensional space
$\range P(A_s,N_t)$. Since $[0,1]$ is compact, there exists a
partition $0=s_0<\ldots<s_n=1$ and $t_k\in[s_k,s_{k+1}]$,
$k=0,\ldots,n-1$ such that
$[s_k,s_{k+1}]\subset(t_k-\delta(t_k),t_k+\delta(t_k))$ for each
$k=0,\ldots,n-1$.

\begin{definition}\label{d:sf-sa} The {\em spectral flow} $\SF\{A_s\}$
of the family $\{A_s\}$, $0\le s\le 1$ is defined by
\begin{equation}\label{e:sp-flow-sa}\SF\{A_s\}
:=\sum_{k=0}^{n-1}
\Bigl(m^-\bigl(A(s_k,t_k)\bigr) - m^-\bigl(A(s_{k+1},t_k)\bigr)\Bigr).
\end{equation}
\end{definition}

\medskip

After carefully examining the above definition, inspired by
\cite{Ph96}, we find that the necessary data for defining any
spectral flow are the following:
\begin{itemize}
\item a co-oriented bounded real $1$-dimensional regular $C^1$
submanifold $\ell$ of $\C$ without boundary (we call such an
$\ell$ {\em admissible}, and denote by $\ell\in\Aa(\C)$); \item a
Banach space $X$; \item and a continuous family of admissible
operators $A_s$, $0\le s\le 1$ in $\Aa_{\ell}(X)$.
\end{itemize}
Here we define $A\in\Cc(X)$ to be {\em admissible} with respect to
$\ell$, if there exists a bounded open neighborhood $N$ of $\ell$
in $\C$ with $C^1$ boundary $\partial N$ such that (i) $\partial
N\cap\sigma(A)=\emptyset$; (ii) $N\cap\sigma(A)\subset \ell$ is a
finite set; and (iii) $P^0_{\ell}(A):= P(A, N)$ is a
finite rank projection.

Note that $P^0_{\ell}(A)$ does not depend on the specific choice of $N$.
We call $\nu_{h,\ell}(A):=\dim\range P^0_{\ell}(A)$ the {\em
hyperbolic nullity} of $A$ with respect to $\ell$. We denote by
$\Aa_{\ell}(X)$ the set of closed admissible operators with
respect to $\ell$. It is an open subset of $\Cc(X)$.

Similarly as before, we can define the spectral flow
$\SF_{\ell}\{A_s\}$. It counts the number of spectral lines of
$A_s$ coming from the negative side of $\ell$ to the non-negative
side of $\ell$. For the details, see \cite{ZhLo99}.

\begin{examples}\label{ex:co-orented-curve}
a) In the above self-adjoint case,
$\ell=\sqrt{-1}(-\epsilon,\epsilon)$ ($\epsilon>0$) with
co-orientation from left to right. Then a self-adjoint operator
$A$ is admissible with respect to $\ell$ if and only if $A$ is
Fredholm.

\noi b) Another important case is that
$\ell=(1-\epsilon,1+\epsilon)$ ($\epsilon\in(0,1)$) with
co-orientation from downward to upward, and all $A_s$ unitary. A
unitary operator $A$ is admissible with respect to $\ell$ if and
only if $A-I$ is Fredholm.
\end{examples}

The spectral flow has the following properties (cf. \cite{Ph96}
and Lemma 2.6 and Proposition 2.2 in \cite{ZhLo99}).

\begin{proposition}\label{p:sp-flow} Let $\ell\in\Aa(\C)$ be admissible
and let $\{A_s\}$, $0\le s\le 1$ be a curve in $\Aa_{\ell}(X)$.
Then the spectral flow $\SF_{\ell}\{A_s\}$ is well defined, and
the following holds:
\begin{enumerate}
\item {\bf Catenation.} Assume $t\in [0,1]$. Then we have
\begin{equation}\label{e:sf-catenation}
\SF_{\ell}\{A_s;0\le s\le t\}+\SF_{\ell}\{A_s;t\le s\le
1\}=\SF_{\ell}\{A_s;0\le s\le 1\}.
\end{equation}

\item {\bf Homotopy invariance.} Let $A(s,t)$,
$(s,t)\in[0,1]\times[0,1]$ be a continuous family in
$\Aa_{\ell}(X)$. Then we have
\begin{equation} \label{e:sf-homotopy}
\SF_{\ell}\{A(s,t);(s,t)\in\partial([0,1]\times[0,1])\}=0.\end{equation}

\item {\bf Endpoint dependence for Riesz continuity.} Let
$\Bb^{\sa}(X)$, respectively $\Cc^{\sa}(X)$ denote the spaces of
bounded, respectively closed self-adjoint operators in $X$. Let
\[
\begin{matrix}
R&:&\Cc^{\sa}&\to&\Bb^{\sa}(X)\\
\ &\, & A&\mapsto & A(A^2+I)^{-\frac{1}{2}}
\end{matrix}
\]
denote the {\em Riesz transformation}. Let $A_s\in\Cc^{\sa}(X)$
for $s\in [0,1]$. Assume that $\{R(A_s)\}$, $0\le s\le 1$ is a
continuous family. If $m^-(A_0)<+\infty$, then $m^-(A_1)<+\infty$
and we have
\begin{equation}\label{sf-general}
\SF\{A_s\}=m^-(A_0)-m^-(A_1).
\end{equation}

\item {\bf Product.} Let $\{P_s\}$ be a curve of projections on
$X$ such that $P_sA_s\subset A_sP_s$ for all $s\in[0,1]$. Set
$Q_s=I-P_s$. Then we have $P_sA_sP_s\in \Aa_{\ell}(\range P_s)\<
\Cc(\range P_s)$, $Q_sA_sQ_s\in \Aa_{\ell}(\range Q_s)\<
\Cc(\range Q_s)$, and
\begin{equation} \label{sf-product}
\SF_{\ell}\{A_s\}=\SF_{\ell}\{P_sA_sP_s\}+\SF_{\ell}\{Q_sA_sQ_s\}.
\end{equation}

\item {\bf Bound.} For $A\in\Aa_{\ell}(X)$, there exists a
neighborhood $\Nn$ of $A$ in $\Cc(X)$ such that
$\Nn\subset\Aa_{\ell}(X)$, and for curves $\{A_s\}$ in $\Nn$ with
endpoints $A_0=:A$ and $A_1=:B$, the {\em relative Morse index}
$I_{\ell}(A,B):= -\SF_{\ell}\{A_s,0 \, ; \, \le s\le 1\}$ is well
defined and satisfies
\begin{equation}\label{sf-bound}
0\le I_{\ell}(A,B)\le\nu_{h,\ell}(A)-\nu_{h,\ell}(B).
\end{equation}

\item {\bf Reverse orientation.} Let $\hat {\ell}$ denote the
curve $\ell$ with opposite co-orientation. Then we have
\begin{equation}\label{e:sf-reverse}
\SF_{\ell}\{A_s\}+\SF_{\hat
{\ell}}\{A_s\}=\nu_{h,\ell}(A_1)-\nu_{h,\ell}(A_0).
\end{equation}

\item {\bf Zero.} Suppose that $\nu_{h,\ell}(A_s)$ is constant for
$s\in [0,1]$. Then $\SF_{\ell}\{A_s\}=0$.

\item {\bf Invariance.} Let $\{T_s\}_{s\in[0,1]}$ be a curve of
bounded invertible operators. Then we have
\begin{equation}\label{sf-invariance}\SF_{\ell}\{T_s^{-1}A_sT_s\}=\SF_{\ell}\{A_s\}.
\end{equation}

\end{enumerate}
\end{proposition}

\smallskip

Now we give a method of calculating the spectral flow of
differentiable curves, inspired among others by J.J. Duistermaat
\cite{Du76} and J. Robbin and D. Salamon \cite{RoSa93}.

\begin{definition}\label{d:crossing} Let $\ell\in\Aa(\C)$ be admissible and
$\{A_s\}_{s\in[0,1]}$ be a curve in $\Aa_{\ell}(X)$.
\begin{enumerate}
\item A {\em crossing} for $A_s$ is a number $t\in [0,1]$ such
that $\nu_{h,\ell}(A_t)\ne 0$.

\item Set $P_s=P^0_{\ell}{A_s}$. A crossing $t$ is called {\em
regular} if $\dom(A_s)=D$ fixed for $s$ near $t$, $A_sx$ is
differentiable at $s=t$ for all $x\in D$, and $P_t\dot{A_t}P_t$ is
{\em hyperbolic}, i.e. $\nu_{h,\ell}(P_t\dot{A_t}P_t)=0$, where
$\dot{A_s} $ is the unbounded operator with domain $D$ defined by
$$\dot{A_s}x=\frac{d}{ds}A_sx$$
for all $x\in D$.

\item A crossing $t$ is called {\em simple} if it is regular and
$\nu_{h,\ell}(A_t)=1$.

\end{enumerate}\end{definition}

\begin{proposition}[cf. Theorem  4.1 of \cite{ZhLo99}]\label{p:cal-sf}
Let $X$ be a Banach space and $\ell=\sqrt{-1}(-\epsilon,\epsilon)$
($\epsilon>0$) with co-orientation from left to right. Let $A_s$,
$-\epsilon\le s\le \epsilon$ ($\epsilon>0$), be a curve in
$\Aa_{\ell}(X)$. Suppose that $0$ is a regular crossing of $A_s$.
Set $P=P^0_{\ell}(A_0)$, $A=A_0$ and $B=\dot A_s|_{s=0}$. Assume
that
\begin{equation}\label{e:cal-sf-com}
P(AB-BA)P=0.
\end{equation}
Then there is a $\delta\in(0,\epsilon)$ such that
$\nu_{h,\ell}(A_s)=0$ for all $s\in[-\delta,0)\cup (0,\delta]$ and
\begin{eqnarray}
& &\label{e:cal-sf1}\SF_{\ell}\{A_s ; 0\le s\le\delta\}=-m^{-}(PBP),\\
& &\label{e:cal-sf2}\SF_{\ell}\{A_s ; -\delta\le s\le
0\}=m^+(PBP).
\end{eqnarray}
Here we denote by $m^+(PBP)$ $($$m^-(PBP)$$)$ the total algebraic
multiplicity of eigenvalues of $PBP$ with positive $($negative$)$
imaginary part respectively.
\end{proposition}

\smallskip

\subsection{Symplectic functional analysis and Maslov index}\label{ss:symplectic}

A main feature of symplectic analysis is the study of the {\em
Maslov index}. It is an intersection index between a path of
Lagrangian subspaces with the {\em Maslov cycle}, or, more
generally, with another path of Lagrangian subspaces. The Maslov
index assigns an integer to each continuous path of Fredholm pairs
of Lagrangian subspaces of a fixed Hilbert space with continuously
varying symplectic structures.

Firstly we define symplectic Hilbert spaces and Lagrangian
subspaces.

\begin{definition}\label{d:algebraic_symplectic_space}
Let $H$ be a complex vector space. A mapping
\[
  \w:H\times H\too \C
\]
is called a (weak) {\em symplectic form} on $H$, if it is
sesquilinear, skew-hermitian, and non-degenerate, i.e.,

\noi (i) $\w(x,y)$ is linear in $x$ and conjugate linear in $y$;

\noi (ii) $\w(y,x)=-\ol{\w(y,x)}$;

\noi (iii) $H^{\w} := \{x\in H \mid \w(x,y)=0\text{ for all $y\in
H$}\} = \{0\}$.

\noi Then we call $(H,\w)$ a {\em complex symplectic vector
space}.
\end{definition}

\begin{definition}\label{d:lagrangian}
Let $(H,\w)$ be a complex symplectic vector space.

\noi (a) The {\em annihilator} of a subspace ${\la}$ of $H$ is
defined by
\[
{\la}^{\w} := \{y\in H \mid \w(x,y)=0 \text{ for all $x\in
{\la}$}\} .
\]

\noi (b) A subspace ${\la}$ is called {\em isotropic}, {\em
co-isotropic}, or {\em Lagrangian} if
\[
{\la} \,\<\, {\la}^{\w}\,,\quad {\la}\,\>\, {\la}^{\w}\,,\quad
{\la}\,=\, {\la}^{\w}\,,
\]
respectively.

\noi (c) The {\em Lagrangian Grassmannian} $\Ll(H,\w)$ consists of
all Lagrangian subspaces of $(H,\w)$.
\end{definition}

\begin{definition}\label{d:sym-hil}Let $H$ be a complex Hilbert space. A
mapping $\omega:H\times H\to \C$ is called a (strong) {\em
symplectic form} on $H$, if $\omega(x,y)=\lla Jx,y\rra_H$ for some
bounded invertible skew-adjoint operator $J$. $(H,\omega)$ is
called a (strong) {\em symplectic Hilbert space}.
\end{definition}

Before giving a rigorous definition of the Maslov index, we fix
the terminology and give a simple lemma.

We recall:

\begin{definition}\label{d:fredholm_pair}
(a) The space of (algebraic) {\em Fredholm pairs} of linear
subspaces of a vector space $H$ is defined by
\begin{equation}\label{e:fp_alg}
\Ff^{2}_{\operatorname{alg}}(H):=\{({\lambda},{\mu})\mid
\dim\left( {\lambda}\cap{\mu}\right)  <+\infty \text{ and } \dim
\bigl(H/(\lambda+\mu)\bigr)<+\infty\}
\end{equation}
with
\begin{equation}\label{e:fp_index}
\index(\la,\mu):=\dim(\la\cap\mu) - \dim(H/(\la+\mu)).
\end{equation}

\noi (b) In a Banach space $H$, the space of (topological) {\em
Fredholm pairs} is defined by
\begin{multline}\label{e:fp}
\Ff^{2}(H):=\{({\lambda},{\mu})\in\Ff^2_{\operatorname{alg}}(H)\mid
{\lambda},{\mu}, \tand {\lambda}+{\mu} \subset H \text{ closed}\}.
\end{multline}
\end{definition}

We need the following well-known lemma (see, e.g., \cite[Lemma
1.7]{BoZh04}).

\begin{lemma}\label{l:strong_symplectic}
Let $(H,\w)$ be a (strong) symplectic Hilbert space. Then
\begin{enumerate}
\item there is a $1$-$1$ correspondence between the space
$$\Uu^J=\{U\in\Bb\bigl(H^+,H^-\bigr)|\,\,U^*J|_{H^-}U=-J|_{H^+}\}$$
and $\Ll(H,\w)$ under the mapping $U\to L:= \GGG(U)$ $($= graph of
$U$$)$, where $H^{\pm}=H^{\mp}(\sqrt{-1}J)$ in the sense of the
decomposition \eqref{e:spectral-decomposition};

\item if $U,V\in\Uu^J$ and $\lambda:=\GGG(U)$, $\mu:=\GGG(V)$,
then $(\lambda,\mu)$ is a Fredholm pair if and only if $U-V$, or,
equivalently, $UV\ii-I$ is Fredholm. Moreover, we have a natural
isomorphism
\begin{equation}\label{e:unitary_counting}
\ker(UV\ii-I) \simeq \lambda\cap \mu\,.
\end{equation}
\end{enumerate}

\end{lemma}

\begin{definition}\label{d:maslov} Let $(H,\lla\cdot,\cdot\rra_s)$,
$s\in[0,1]$ be a continuous family of Hilbert spaces, and
$\omega_s(x,y)=\lla J_sx,y\rra_s$ be a continuous family of
symplectic forms on $H$, i.e., $\{A_{s,0}\}$ and $\{J_s\}$ are two
continuous families of bounded invertible operators, where
$A_{s,0}$ is defined by
$$\lla x,y\rra_s=\lla A_{s,0}x,y\rra_0\quad\text{for all}\;x,y\in H.$$
Let $\{(\lambda_s,\mu_s)\}$ be a continuous family of Fredholm
pairs of Lagrangian subspaces of
$(H,\lla\cdot,\cdot\rra_s,\omega_s)$. Then there is a continuous
splitting
\begin{equation}\label{e:s-split}H=H_s^-(\sqrt{-1}J_s)\oplus
H_s^+(\sqrt{-1}J_s) \end{equation} associated to the self-adjoint
operator $\sqrt{-1}J_s\in\Bb(H,\lla\cdot,\cdot\rra_s)$ for each
$s\in[0,1]$. By Lemma \ref{l:strong_symplectic},
$\lambda_s=\GGG_s(U_s)$ and $\mu_s=\GGG_s(V_s)$ with $U_s$,
$V_s\in \Uu^{J_s}$, where $\GGG_s$ denotes the graph associated to
the splitting (\ref{e:s-split}). We define the {\em Maslov index}
$\Mas\{\lambda_s,\mu_s\}$ by
\begin{equation}\label{e:maslov}\Mas\{\lambda_s,\mu_s\}=-\SF_{\ell}\{U_sV_s^{\ii}\},
\end{equation}
where $\ell:=(1-\epsilon,1+\epsilon)$ with, $\epsilon\in(0,1)$ and
with upward co-orientation.
\end{definition}

\begin{rem}\label{r:arnold}
For finite-dimensional $H$, constant $\mu_s=\mu_0$, and a loop
$\{\la_s\}$, i.e., for $\la_0=\la_1$\,, we notice that
$\Mas\{\la_s,\mu_s\}$ is the winding number of the closed curve
$\{\det(U_s\ii V_0)\}_{s\in [0,1]}$\,. This is the original
definition of the Maslov index as explained in Arnol'd,
\cite{Ar67}.
\end{rem}

\begin{lemma}\label{l:independent of norm}The Maslov index is
independent of the choice of the complete inner product of $H$.
\end{lemma}

\begin{proof} Let $\lla\cdot,\cdot\rra_{s,k}$, $s\in[0,1]$ with $k=0,1$ be
two continuous families of complete inner products of $H$. We
define
$$\lla\cdot,\cdot\rra_{s,t}=(1-t)\lla\cdot,\cdot\rra_{s,0}+t\lla\cdot,\cdot\rra_{s,1}$$
for each $(s,t)\in[0,1]\times[0,1]$. Let $(\lambda_s,\mu_s)$ be a
continuous family of Fredholm pairs of Lagrangian subspaces of
$(H,\omega_s)$. For each inner product
$\lla\cdot,\cdot\rra_{s,t}$, we denote by $U_{s,t}$ and $V_{s,t}$
the associated generated "unitary" operators of $\lambda_s$ and
$\mu_s$ respectively. We also denote by $\Mas_t$ the Maslov index
defined with $\lla\cdot,\cdot\rra_{s,t}$ for each $t\in[0,1]$. By
Proposition \ref{p:sp-flow} we have
\begin{eqnarray*}
\Mas_0\{\lambda_s,\mu_s\}&-&\Mas_1\{\lambda_s,\mu_s\}\\
&=&-\SF_{\ell}\{U_{s,0}V_{s,0}^{\ii}\}+\SF_{\ell}\{U_{s,1}V_{s,1}^{\ii}\}\\
&=&-\SF_{\ell}\{U_{s,t}V_{s,t}^{\ii};(s,t)\in\partial\bigl([0,1]\times[0,1]\bigr)\}\\
&=&0.
\end{eqnarray*}
\end{proof}

\smallskip

Now we give a method of using the crossing form to calculate
Maslov indices (cf. \cite{Du76}, \cite{RoSa93}, \cite[Theorem
2.1]{BoFu98}; for a full proof of the following Proposition see
\cite[Corollary 3.1]{Zh01}).

Let $\la=\{\lambda_s\}_{s\in[0,1]}$ be a $C^1$ curve of Lagrangian
subspaces of $H$. Let $W$ be a fixed Lagrangian complement of
$\lambda_t$. For $v\in\lambda_t$ and $|s-t|$ small, define
$w(s)\in W$ by $v+w(s)\in\lambda_s$. The form
\begin{equation}\label{e:crossing1}
Q(\lambda,t):=
Q(\lambda,W,t)(u,v)=\frac{d}{ds}|_{s=t}\omega(u,w(s)),\quad
\forall u,v\in\lambda_t
\end{equation}
is independent of the choice of $W$. Let
$\{(\lambda_s,\mu_s)\}$, $0\le s\le 1$ be a curve of Fredholm pairs of
Lagrangian subspaces of $H$. For $t\in[0,1]$, the {\em crossing
form} $\Gamma(\lambda,\mu,t)$ is a quadratic form on
$\lambda_t\cap\mu_t$ defined by
\begin{equation}\label{e:crossing2}
\Gamma(\lambda,\mu,t)(u,v)=Q(\lambda,t)(u,v)-Q(\mu,t)(u,v),\quad\forall
u,v\in\lambda_t\cap\mu_t\/.
\end{equation}
A {\em crossing} is a time
$t\in[0,1]$ such that $\lambda_t\cap\mu_t\ne\{0\}$. A crossing
is called {\em regular} if $\Gamma(\lambda,\mu,t)$ is
nondegenerate. It is called {\em simple} if it is regular and
$\lambda_t\cap\mu_t$ is one-dimensional.

\begin{proposition} \label{p:cal-mas} Let $(H,\omega)$ be a symplectic Hilbert space
and $\{(\lambda_s,\mu_s)\}$, $0\le s\le 1$ be a $C^1$ curve of
Fredholm pairs of Lagrangian subspaces of $H$ with only regular
crossings. Then we have
\begin{equation} \label{e:cal-mas}
\Mas\{\lambda,\mu\}=m^+(\Gamma(\lambda,\mu,0))
-m^-(\Gamma(\lambda,\mu,1))+\sum_{0<t<1}\sign(\Gamma(\lambda,\mu,t)).
\end{equation}
\end{proposition}

\medskip

\section{Symplectic analysis of symmetric operators}\label{s:sa-extensions}

\subsection{Local stability of weak inner UCP}\label{ss:stable-ucp}

Let $X$ be a complex Hilbert space and $A\in\Cc(X)$ a linear,
closed, densely defined operator in $X$. We assume that $A$ is
symmetric, i.e., $A^*\>A$ where $A^*$ denotes the adjoint
operator. We denote the domains of $A$ by $D_m$ (the {\it minimal}
domain) and of $A^*$ by $D_{\mmax}$ (the {\it maximal} domain).

\begin{definition}\label{d:ucp}
Let $X$ be a Hilbert space and $A\in \Cc(X)$ with $\dom A=D_m$ and
$A^*\> A$. We shall say that the operator $A$ satisfies the {\em
weak inner Unique Continuation Property (UCP)} if $\ker A=\{0\}$.
\end{definition}

It is well known that weak UCP and weak inner UCP can be
established for a large class of Dirac type operators, see the
first author with Wojciechowski \cite[Chapter 8]{BoWo93}, and the
first author with M. Marcolli and B.-L. Wang \cite{BoMaWa02}.
However, it is not valid for all linear elliptic differential
operators of first order as shown by one of the Pli{\'s}
counter-examples \cite{Pl61}. Moreover, one has various quite
elementary examples of linear and non-linear perturbations which
{\it invalidate} weak inner UCP for Dirac operators. Two such
examples are listed in \cite{BoMaWa02}. In the same paper,
however, it was shown that weak UCP is {\it preserved} under
certain `small' perturbations of Dirac type operators. Here we
show an elementary result, namely the local stability of weak
inner UCP.

\begin{lemma}\label{l:stable-ucp}
Let $X$ be a Hilbert space. Let $A_s\in \Cc(X)$, $0\le
s\le 1$ be a family of symmetric operators with $\dom A_s=D_m$ and
$\dom A_s^*=D_{\mmax}$ independent of $s$. Assume that
$\{A_s^*:D_{\mmax}\to X\}$ is a continuous curve of bounded
operators, where the norm on $D_{\mmax}$ is the graph norm induced
by $A_0^*$. If $A_0$ satisfies weak inner UCP and there exists a
self-adjoint Fredholm extension $A_0^*|_D$ of $A_0$, then for all
$s\ll 1$ the operators $A_s^*$ are surjective and the operators
$A_s$ satisfy weak inner UCP.
\end{lemma}

\begin{proof}
By our assumptions, $\ran A_0^*|_D$ is closed and is of finite
codimension. Since $\ran A_0^*|_D\<\ran A_0^*\< X$, the full range
$\ran A_0^*$ is closed. Since $A_0$ satisfies weak inner UCP,
$\ran A_0^*=X$. Then $A_0^*$ is semi-Fredholm. By Theorem IV.5.17
of Kato \cite{Ka80} we have $\ran A_s^*=X$ for $s\ll 1$. Since $A_s$
are symmetric, $A_s$ satisfy weak inner UCP for $s\ll 1$.
\end{proof}

\smallskip

\subsection{Continuity of the family
$\{A_{s,D_s}\}$}\label{ss:op-continuity}

Let $X$ be a complex Hilbert space, and $M,N\subset X$ be two
closed linear subspaces. Let $P_M,P_N$ be the orthogonal
projections onto $M$, $N$ respectively. Then the distance $d(M,N)$
is defined by $d(M,N)=\|P_M-P_N\|$ and called the {\em gap} between $M$ and $N$.
For any two closed operators
$A,B$ on $X$, we define $d(A,B)$ as the distance between their
graphs.

Let $A\in\Cc(X)$ be a linear, closed, densely defined operator in
$X$. By Footnote 1 (page 198), Theorems IV.1.1 and IV.2.14 in
\cite{Ka80}, it is easy to verify the following

\begin{lemma}\label{l:stable} Let $B\in\Bb(\dom(A),X)$ be a bounded operator,
where the norm on $\dom(A)$ is the graph norm $\Gg_A$ induced by
$A$. Let $d:=\|B-A\|_{\Gg_A}<\frac{1}{2}$. Then we have
\begin{enumerate}
\item $B\in\Cc(X)$, and it holds that
$$
(1-2d)\lla x,x\rra_{\Gg_A}\le\lla x,x\rra_{\Gg_B}\le (1+d)^2\lla
x,x\rra_{\Gg_A} \;\text{ for $x\in D$}.
$$

\item $d(B,A)\le \frac{\sqrt{2}d}{(1-d)^{-1}}$.
\end{enumerate}
\end{lemma}

\begin{lemma}\label{l:sub-quotion}Let $X$ be a Hilbert space, and $Y$
be a closed linear subspace of $H$. Then there exists a bijection
between the space of closed linear subspaces of $X$ containing $Y$
and that of closed linear subspaces of $X/Y$ which preserves the
metric.
\end{lemma}

\begin{proof} We view $X/Y$ as $Y^{\bot}$. Let $M,N\subset Y^{\bot}$ be
two closed subspaces and $P_M,P_N$ be the orthogonal projections
onto $M$, $N$ respectively. Then we have
$$d(M+Y,N+Y)=\|P_{M+Y}-P_{N+Y}\|=\|P_M-P_N\|=d(M,N).$$
\end{proof}

From the definition of the gap norm and by some computations we have

\begin{lemma}\label{l:op-continuity}
Let $D_m\subset D_{\mmax}\subset X$ be three Hilbert spaces such
that $D_m$ is a closed subspace of $D_{\mmax}$ and a dense
subspace of $X$. Let $\bigl\{A_s\in\Cc(X)\bigr\}_{s\in[0,1]}$ be a
family of densely defined symmetric operators with domain $D_m$,
and $\bigl\{D_s\bigr\}_{s\in[0,1]}$ be a family of closed
subspaces of $D_{\mmax}$ containing $D_m$. We assume that
$\dom(A_s^*)=D_{\mmax}$, each graph norm $\Gg_s$ of $D_{\mmax}$
induced by $A_s^*$ is equivalent to the original norm $\Gg$ of
$D_{\mmax}$\,, and $\bigl\{A_s^*\in\Bb(D_{\mmax},X)\bigr\}$,
$\bigl\{D_s/D_m\subset D_{\mmax}/D_m\bigr\}$ are two continuous
families. Then $\bigl\{A_{s,D_s}\in\Cc(X)\bigr\}_{s\in[0,1]}$ is a
continuous family of closed operators.
\end{lemma}

\smallskip

\subsection{Continuity of natural Cauchy data spaces}\label{ss:cd-continuity}

In this subsection we generalize the proof of the continuity of
Cauchy data spaces given in \cite[Section 3.3]{BoFu98}.
We need the following

\begin{proposition}[Proposition 3.5 of \cite{BoFu98}]
\label{p:cauchy} Let $X$ be a Hilbert space, and $A\in\Cc(X)$ be a
symmetric operator. Set $D_m=\dom(A)$ and $D_{\mmax}=\dom(A^*)$.
If $A$ admits a self-adjoint Fredholm extension with domain $D$\,,
then the quotient space $D/D_m$ and the natural Cauchy data space
$(\ker A^*+D_m)/D_m$ form a Fredholm pair of Lagrangian subspaces
of the $($strong$)$ symplectic Hilbert space $D_{\mmax}/D_m$
$($introduced above in Subsection $\ref{ss:statement-fixd}$, Item
$(ii))$.
\end{proposition}

\begin{rem}
At present (March 2005), it is not known whether all linear formally
self-adjoint elliptic differential operators of first order over
a compact smooth Riemannian manifold with smooth boundary admit a self-adjoint
Fredholm extension. Recently in \cite{BoLeZh05}, however, that crucial property
has been established under the additional assumption of self-adjoint principal
symbol of the ``tangential operator" at the boundary.
\end{rem}


Now we can prove

\begin{proposition}\label{p:cd_continuity}
Let $X$ be a Hilbert space, and  $D_m\subset D_{\mmax}$ be two
dense linear subspaces of $X$. Let $\{A_s:D_m\to X\}_{s\in [0,1]}$
be a family of closed symmetric densely defined operators in $X$.
We assume that
\begin{enumerate}
\item each $A_s$ admits a self-adjoint Fredholm extension with
domain $D_s$\,;

\item $\dom(A_s^*)=D_{\mmax}$ independent of $s$ and that all graph
norms $\Gg_s$ of $D_{\mmax}$ induced by $A_s^*$ are mutually
equivalent;

\item each $A_s$ satisfies weak inner UCP relative to $D_m$\,; and

\item $\{A_s^*:D_{\mmax}\to X\}$ forms a continuous family of bounded
operators, where the norm on $D_{\mmax}$ is the graph norm $\Gg$
induced by $A_0$.
\end{enumerate}
Then the natural Cauchy data spaces $\bigl(D_m+\ker
A_s^*\bigr)/D_m$ are continuously varying in $D_{\mmax}/D_m$\,.
\end{proposition}

\begin{proof} We denote the
projection of $D_{\mmax}$ onto $D_{\mmax}/D_m$ by $\g$\,. Note that $\ker
A_s^*$ is closed in $D_{\mmax}$\,.

To prove the continuity, we need only to consider the local
situation at $s=0$. First we show that $\{\ker A_s^*\}_{s\in
[0,1]}$ is a continuous family of subspaces of $D_{\mmax}$; then
we show that $\g(\ker A_s^*)$ is a continuous family in
$D_{\mmax}/D_m$.

\smallskip

We consider the bounded operator
\[
\begin{matrix}
F_s: & D_{\mmax}&\too & X \oplus \ker
A_0^*        \\
\ & x& \mapsto & \left(A^*_s(x), P_0x\right)
\end{matrix}\qquad ,
\]
where $P_0:D_{\mmax}\to \ker A_0^*$ denotes the orthogonal
projection of the Hilbert space $D_{\mmax}$ onto the closed
subspace $\ker A_0^*$\,. By definition, the family $\{F_s\}$ is a
continuous family of bounded operators.

Clearly, $F_0$ is injective. Since $\ran A_0^*|_{D_0}\subset \ran
A_0^*\subset X$ and $A_0^*|_{D_0}$ is Fredholm, $\ran A_0^*$ is
closed. From weak inner UCP we get $\ran A_0^*=X$\/. So the
operator $F_0$ is also surjective. This proves that $F_0$ is
invertible with bounded inverse. Then all operators $F_s$ are
invertible for small $s \geq 0$, since $F_s$ is a continuous
family of operators.

Note that
\[
F_s(\ker A_s^*)\subset \{0\}\oplus \ker A_0^*,\quad (
F_s)^{-1}(\{0\}\oplus \ker A_0^*)\subset\ker A_s^*.
\]
Since $F_s$ are invertible for small $s \geq 0$, we have
\begin{equation}\label{e:f-range}
F_s(\ker A_s^*)= \{0\}\oplus \ker A_0^*.
\end{equation}

We define
\[
\f_s:= F_s\ii\circ F_0 :D_{\mmax}\ \cong\ D_{\mmax} \tand \f_s\ii=
F_0\ii\circ F_s :D_{\mmax}\ \cong\ D_{\mmax}
\]
for $s$ small. Since $F_s$ are invertible for small $s \geq 0$,
from \eqref{e:f-range} we obtain that
\begin{equation}\label{e:phi}
\f_s(\ker A_0^*)=\ker A_s^*\,.
\end{equation}

From \eqref{e:phi} we get that
\[
\{P_s:=\f_s P_0\f_s\ii:D_{\mmax}\too \ker A_s^*\}
\]
is a continuous family of projections onto the solution spaces
$\ker A_s^*$. The projections are not necessarily orthogonal, but
can be orthogonalized and remain continuous in $s$ like in
\cite[Lemma 12.8]{BoWo93}. This proves the continuity of the
family $\{\ker A_s^*\}$ in $D_{\mmax}$\,.

Now we must show that $\{\g(\ker A_s^*)\}$ is a continuous family
in the quotient space $D_{\mmax}/D_m$. This is not proved by the formula $\g(\ker
A_s^*)=\g(\f_s(\ker A_0^*))$ alone. We must modify the
endomorphism $\f_s$ of $D_{\mmax}$ in such a way that it keeps the
subspace $D_m$ invariant.

By Proposition (\ref{p:cauchy}), the Cauchy data space
$\g(D_m+\ker A_0^*)$ is closed in $D_{\mmax}/D_m$\/. So $D_m+\ker
A_0^*$ is closed in $D_{\mmax}$. We define a continuous family of
mappings by
\[
\begin{matrix}
\psi_s: D_{\mmax} = &D_m  +  \ker A_0^* & \ + \ & (D_m+\ker
A_0^*)^\perp & \
\too \ & D_{\mmax}\\
\ & x+y & + & z&\mapsto&x+\f_s(y)+z
\end{matrix}
\]
with $\psi_0=\id$. Hence all $\psi_s$ are invertible for $s\ll 1$,
and $\psi_s(D_m)=D_m$ for such small $s$. Hence we obtain a
continuous family of mappings $\{\wt\psi_s:D_{\mmax}/D_m\to
D_{\mmax}/D_m\}$ with $\wt\psi_s(\g(\ker A_0^*))=\g(\ker A_s^*)$.
From that we obtain a continuous family of projections as above.
\end{proof}

\begin{rem}
From the preceding arguments it also follows that the Cauchy data
spaces form a differentiable family, if $\{A_s^*\}$ is a
differentiable family.
\end{rem}

\smallskip

\subsection{Proof of the spectral flow formula}\label{s:proof}

We begin with a simple case.

\begin{lemma}\label{l:bofu98}
Let $X$ be a Hilbert space, and $A\in\Cc(X)$ be a symmetric
operator with $\dom(A)=D_m$ and $\dom(A^*)=D_{\mmax}$. Let
$A_D:=A^*|_{D}$ be a self-adjoint Fredholm extension of $A$. We
assume that $A$ satisfies weak inner UCP. Then there exists an
$\epsilon>0$ such that $A_D+aI$ is Fredholm and satisfies weak
inner UCP for each $a\in[0,\epsilon]$. Let $\g:D_{\mmax}\to
D_{\mmax}/D_m$ be the natural projection. Then we have
\[
\SF\{A_D+aI; a\in [0,\epsilon]\} = -\Mas \{\g(D),\g(\ker
(A^*+aI));a\in [0,\epsilon]\}.
\]
\end{lemma}

\begin{proof} By the definition of the spectral flow we have
\begin{equation}\label{e:sf-mas1}
\SF\{A_D+aI; a\in [0,\epsilon]\}=\sum_{a\in(0,\epsilon]}\dim\ker
(A_D+aI).
\end{equation}

Let $\omega$ be the Green form on $D_{\mmax}$ induced by $A^*$.
Let $W\in\Ll(D_{\mmax}/D_m)$ be a Lagrangian complement of
$\g(\ker(A^*+a_0I))$. By Proposition \ref{p:cd_continuity},
$\g(\ker(A^*+aI))$ and $\ker(A^*+aI)$ are two differentiable
families. For each $y(a_0)\in\ker(A_D+a_0I)$, there exists a
continuous family $w(a)\in W+D_m$, $|a-a_0|$ small, such that
$w(a_0)=0$ and $y(a):=y(a_0)+w(a)\in \ker(A^*+aI)$. Since
$A^*(y(a))=-ay(a)$ and the family $\{y(a)\}$ is continuous in
$D_{\mmax}$, the family $\{y(a)\}$ is also continuous in $X$. For
all $x(a_0)\in\ker(A_D+a_0I)$, we have
\begin{eqnarray*}
& &\omega(\g(x(a_0)),\g(w(a))\\
&=&\lla A^*(x(a_0)),y(a)-y(a_0)\rra-\lla x(a_0),A^*(w(a))\rra \\
&=&\lla -a_0x(a_0),y(a)-y(a_0)\rra-\lla
x(a_0),A^*(y(a))-A^*(y(a_0))\rra\\
&=&\lla -a_0x(a_0),y(a)-y(a_0)\rra-\lla
x_(a_0),-ay(a)+a_0y(a_0)\rra\\
&=&(a-a_0)\lla x(a_0),y(a)\rra
\end{eqnarray*}
Let the crossing forms $Q$ and $\Gamma$ be defined by
(\ref{e:crossing1}) and (\ref{e:crossing2}) respectively. Then we
have $Q(\g(\ker (A^*+aI)),a_0)(\g(x(a_0)),\g(y(a_0))=\lla
x(a_0),y(a_0)\rra$ and
\[
\Gamma(\g(D),\g(\ker (A^*+aI)),a_0)(\g(x(a_0)),\g(y(a_0))=-\lla
x(a_0),y(a_0)\rra.
\]
By Proposition \ref{p:cal-mas} we have
\begin{multline}\label{e:sf-mas2}
\Mas \{\g(D),\g(\ker (A^*+aI));a\in
[0,\epsilon]\}\\
=-\sum_{a\in(0,\epsilon]}\dim\ker \bigl(A_D+aI\bigr).
\end{multline}
Combine equations (\ref{e:sf-mas1}), (\ref{e:sf-mas2}), and our
lemma follows.
\end{proof}

Now our main result follows at once.

\begin{proof}[Proof of Theorem \ref{t:fixm}]
By Lemma \ref{l:stable-ucp}, for each $s_0$ there exists an
$\epsilon(s_0)>0$ such that the operators $A_s+aI$ satisfy weak
inner UCP for all $s,a$ with $\abs{s-s_0}, \abs{a}<\epsilon(s_0)$.
Here we use the continuity of the family $\bigl\{A_s^*\}$ as
bounded operators from $D_{\mmax}$ to $X$. Since $[0,1]$ is
compact and $A_{s,D_s}$ are Fredholm operators for all
$s\in[0,1]$, there exists an $\epsilon>0$ such that the operators
$A_s+aI$ satisfy weak inner UCP and $A_{s,D_s}+aI$ are Fredholm
operators for all $s\in [0,1]$ and $\abs{a}<\epsilon$.

We only need to prove the formula \eqref{e:fixd} in a small
interval $[s_0,s_1]$. We consider the two-parameter families
\[
\{A_{s,D_s}+aI\} \tand \{\g(D_s),\g(\ker A_s^*+aI)\}
\]
for $s\in[s_0,s_1]$ and $a\in[0,\epsilon]$. Because of the
homotopy invariance of spectral flow and Maslov index, both
integers must vanish for the boundary loop going counter clockwise
around the rectangular domain from the corner point $(s_0,0)$ via
the corner points $(s_1,0)$, $(s_1,\epsilon)$, and
$(s_0,\epsilon)$ back to $(s_0,0)$.

Moreover, for $s_1$ sufficiently close to $s_0$ we can choose
$\epsilon$ sufficiently small so that $\ker(A_{s,D_s}+\epsilon
I)=\{0\}$ for all $s\in[s_0,s_1]$. Hence, spectral flow and Maslov
index must vanish on the top segment of our box.

Finally, by the preceding lemma, the left and the right side
segments of our curves yield vanishing sum of spectral flow and
Maslov index. So, by additivity under catenation, our assertion
follows.
\end{proof}

\bigskip
\addtocontents{toc}{\medskip\noi}

\end{document}